# RANDOM SYSTEMS WITH COMPLETE CONNECTIONS AND THE GAUSS PROBLEM FOR THE REGULAR CONTINUED FRACTIONS

Dan Lascu, Ion Coltescu

Naval Academy "Mircea cel Batran", Constanta

*lascudan@gmail.com, icoltescu@yahoo.com*

**ABSTRACT**

*This paper present the important role that random system with complete connections played in solving the Gauss problem associated to the regular continued fractions. Hence, using the ergodic behavior of homogeneous random system with complete connections, we will solve a Gauss – Kuzmin type theorem.*

## 1. DEFINITIONS

Let $\Omega$ denote the collection of irrational numbers in the unit interval $I = [0,1]$. Write $x \in \Omega$ as a regular continued fraction (RCF)

$$x = \cfrac{1}{a_1(x) + \cfrac{1}{a_2(x) + \ddots}} := [a_1(x), a_2(x), \ldots], \qquad (1.1)$$

where the $\mathbf{N}_+$-valued functions on $\Omega$, $a_k(x)$, are unique and called incomplete quotients of $x \in \Omega$. We will usually drop the dependence on $x$ in the notation of incomplete quotients $a_n$, $n \in \mathbf{N}_+ := \{1, 2, \ldots\}$.

Define on $I$ the transformation $\tau$ as follows

$$\tau(x) := \tau([a_1, a_2, \ldots]) = [a_2, a_3, \ldots], \ x \neq 0; \ \tau(0) := 0. \qquad (1.2)$$

It follows from (1.1) and (1.2) that for $x \neq 0$ we have

$$x = \frac{1}{a_1 + \tau(x)}. \qquad (1.3)$$

Consequently, we can write the transformation $\tau$ of $[0,1)$ as

$$\tau(x) := \frac{1}{x} - \left\lfloor \frac{1}{x} \right\rfloor, \ x \neq 0, \qquad (1.4)$$

where $\lfloor \cdot \rfloor$ denotes the floor function. Usually, this transformation is called the *Gauss map*. For $x \neq 0$, we get

$$a_1 = \left\lfloor \frac{1}{x} \right\rfloor \text{ and } a_n = a_1(\tau^{n-1}(x)), \ n \in \mathbf{N}_+, \ n \geq 2, \qquad (1.5)$$

with $\tau^0(x) = x$ and $\tau^n(x) = \tau(\tau^{n-1}(x))$. For any $n \in \mathbf{N}_+$, writing

$$[x_1] = \frac{1}{x_1}, \ [x_1, x_2, \ldots, x_n] = \frac{1}{x_1 + [x_2, \ldots, x_n]}, \ n \geq 2$$

for arbitrary indeterminates $x_i$, $1 \leq i \leq n$, we have

$$x = \lim_{n \to \infty} [a_1, \ldots, a_n], \ x \in \Omega. \qquad (1.6)$$

The probability structure of the sequence $(a_n)_{n \in \mathbf{N}_+}$ under $\lambda$, the Lebesgue measure, is described by the equations

$$\lambda(a_1 = i) = \frac{1}{i(i+1)}, \ i \in \mathbf{N}_+, \quad (1.7)$$

$$\lambda(a_{n+1} = i | a_1, \ldots, a_n) := P_i(s_n), \ i, n \in \mathbf{N}_+, \quad (1.8)$$

where

$$s_n = [a_n, \ldots, a_1] \quad (1.9)$$

and

$$P_i(x) = \frac{x+1}{(x+i)(x+i+1)}, \ i \in \mathbf{N}_+, \ x \in I. \quad (1.10)$$

Thus, under $\lambda$, the sequence $(a_n)_{n \in \mathbf{N}_+}$ is neither independent nor Markovian. Then, if $\mathcal{B}_I$ denotes the $\sigma$-algebra of Borel subsets of $I$, there is a probability measure $\gamma$ on $\mathcal{B}_I$

$$\gamma(A) = \frac{1}{\log 2} \int_A \frac{dx}{x+1}, \ A \in \mathcal{B}_I, \quad (1.11)$$

called *Gauss' measure*, which makes $(a_n)_{n \in \mathbf{N}_+}$ into a strictly stationary sequence. Moreover, $\gamma$ is $\tau$-invariant, i.e.

$$\gamma(\tau^{-1}(A)) = \gamma(A), \ A \in \mathcal{B}_I. \quad (1.12)$$

## 2. RANDOM SYSTEM WITH COMPLETE CONNECTIONS

**Definition 2.1** A quadruple $\{(W, \mathcal{W}), (X, \mathcal{X}), u, P\}$ is named a homogeneous random system with complete connections (RSCC) if
  (i) $(W, \mathcal{W})$ and $(X, \mathcal{X})$ are arbitrary measurable spaces;
  (ii) $u : W \times X \to W$ is a $(\mathcal{W} \otimes \mathcal{X}, \mathcal{W})$ - measurable function;
  (iii) $P$ is a transition probability function from $(W, \mathcal{W})$ to $(X, \mathcal{X})$.

Next, denote the element $(x_1, x_2, \ldots, x_n) \in X^n$ with $x^{(n)}$.

**Definition 2.2** The functions $u^{(n)} : W \times X^n \to W$, $n \in \mathbf{N}_+$ are defined as follows

$$u^{(n+1)}(w, x^{(n+1)}) = \begin{cases} u(w, x), & \text{if } n = 0 \\ u(u^{(n)}(w, x^{(n)}), x_{n+1}), & \text{if } n \in \mathbf{N}_+. \end{cases} \quad (2.1)$$

By convention, we will write $wx^{(n)}$ instead of $u^{(n)}(w, x^{(n)})$.

**Definition 2.3** The transition probability functions $P_r$, $r \in \mathbf{N}_+$ are defined by

$$P_r(w, A) = \begin{cases} P(w, A), & \text{if } r = 1, \\ \sum_{x_1 \in X} P(w, x_1) \sum_{x_2 \in X} P(wx_1, x_2) \ldots \sum_{x_r \in X} P(wx^{(r-1)}, x_r) \chi_A(x^{(r)}), & \text{if } r > 1, \end{cases} \quad (2.2)$$

for any $w \in W$, $r \in \mathbf{N}_+$ and $A \in \mathcal{X}^r$, where $\chi_A$ is the characteristic function of the set $A$.

**Definition 2.4** Assume that $X^0 \times A = A$. Then we define

$$P_r^n(w, A) = P_{n+r-1}(w, X^{n-1} \times A), \quad (2.3)$$

for any $w \in W$, $n, r \in \mathbf{N}_+$ and $A \in \mathcal{X}^r$.

**Theorem 2.5 (Existence theorem)** Let $\{(W,\mathcal{W}),(X,\mathcal{X}),u,P\}$ be a homogeneous RSCC and let $w_0 \in W$. Then there exist a probability space $(\Omega,\mathcal{K},P_{w_0})$ and two chains of random variables $(\xi_n)_{n\in\mathbf{N}_+}$ and $(\zeta_n)_{n\in\mathbf{N}}$ defined on $\Omega$ with values in $X$ and $W$ respectively, such that

(i) (a) $P_{w_0}\left((\xi_n,\ldots,\xi_{n+r-1}) \in A\right) = P_r^n(w_0, A)$,

(b) $P_{w_0}\left((\xi_{n+m},\ldots,\xi_{n+m+r-1}) \in A \mid \xi^{(n)}\right) = P_r^m\left(w_0\xi^{(n)}, A\right)$ $\quad P_{w_0}$-a.e.

(c) $P_{w_0}\left((\xi_{n+m},\ldots,\xi_{n+m+r-1}) \in A \mid \xi^{(n)},\zeta^{(n)}\right) = P_r^m(\zeta_n, A)$ $\quad P_{w_0}$-a.e.

for any $n,m,r \in \mathbf{N}_+$ and $A \in \mathcal{X}^r$, where $\xi^{(n)}$, $\zeta^{(n)}$ denote the random vectors $(\xi_1,\ldots,\xi_n)$ and $(\zeta_1,\ldots,\zeta_n)$ respectively.

(ii) $(\zeta_n)_{n\in\mathbf{N}}$ is a homogeneous Markov chain with initial distribution concentrated in $w_0$ and with the transition operator $U$ defined by

$$Uf(w) = \sum_{x\in X} P(w,x) f(wx), \qquad (2.4)$$

for any $f$ real $W$-measurable and bounded function.

**Remark 2.6** Letting $m = r = 1$ in (i)(b) we obtain

$$P_{w_0}\left(\xi_{n+1} \in A \mid \xi^{(n)}\right) = P\left(w_0\xi^{(n)}, A\right) \qquad P_{w_0}\text{-a.e.} \qquad (2.5)$$

that is the conditioned distribution of $\xi_{n+1}$ by the past depends actually by this, through $u^{(n)}$. This fact justifies the name of *chain of infinite order* or *chain with complete connections* used for $(\xi_{n+1})_{n\in\mathbf{N}_+}$.

**Remark 2.7** On account of (2.4) we have

$$U^n f(w) = \sum_{x^{(n)} \in X^n} P_n\left(w, x^{(n)}\right) f\left(wx^{(n)}\right), \quad n \in \mathbf{N}_+ \qquad (2.6)$$

for any $f$ real $W$-measurable and bounded function.

**Remark 2.8** The transition probability function of the Markov chain $(\zeta_n)_{n\in\mathbf{N}_+}$ is

$$Q(w,A) = \sum_{x\in X} P(w,x)\chi_A(wx) = P(w, A_w), \qquad (2.7)$$

where $A_w = \{x \in X : wx \in A\}$, $w \in W$. It follows that the transition probability after $n$ paths of the Markov chain $(\zeta_n)_{n\in\mathbf{N}}$ is

$$Q^n(w, A) = P_n\left(w, A_w^{(n)}\right), \qquad (2.8)$$

where $A_w^{(n)} = \left\{x^{(n)} \in X : wx^{(n)} \in A\right\}$.

**Definition 2.9** Let $Q_n$ be the transition probability function defined by

$$Q_n(w,A) = \frac{1}{n}\sum_{k=1}^n Q^k(w,A) \qquad (2.9)$$

for any $w \in W$ and $A \in \mathcal{W}$.

**Definition 2.10** Let $U_n$ be the Markov operator associated with $Q_n$. Then

(i) If there exists a linear bounded operator $U^\infty$ from $L(W)$ to $L(W)$ such that

$$\lim_{n\to\infty} \|U_n f - U^\infty f\| = 0, \qquad (2.10)$$

for any $f \in L(W)$ with $\|f\| = 1$, we say that $U$ is *ordered*.

(ii) If

$$\lim_{n\to\infty} \|U^n f - U^\infty f\| = 0, \qquad (2.11)$$

for any $f \in L(W)$ with $\|f\| = 1$, we say that $U$ is *aperiodic*.

(iii) If $U$ is ordered and $U^{\infty}(L(W))$ is one-dimensional space, then $U$ is named *ergodic* with respect to $L(W)$.

(iv) If $U$ is ergodic and aperiodic, then $U$ is named *regular* with respect to $L(W)$ and the corresponding Markov chain has the same name.

**Definition 2.12** If $\{(W,\mathcal{W}),(X,\mathcal{X}),u,P\}$ is a RSCC which satisfies the properties

(i) $(W,d)$ is a metric separable space;

(ii) $r_1 < \infty$, where

$$r_k = \sup_{w' \neq w''} \sum_{X^k} P_k(w, x^{(k)}) \frac{d(w'x^{(k)}, w''x^{(k)})}{d(w', w'')}, \ k \in \mathbf{N}_+; \tag{2.12}$$

(iii) $r_1 < \infty$, where

$$R_1 = \sup_{A \in \mathcal{X}} \sup_{w' \neq w''} \frac{|P(w', A) - P(w'', A)|}{d(w', w'')}; \tag{2.13}$$

(iv) there exists $k \in \mathbf{N}_+$ such that $r_k < 1$,

where $d(x,y) = |x - y|$, for any $x, y \in I$, then we say that this RSCC is a *RSCC with contraction*.

**Theorem 2.13** Let $(W,d)$ be a compact space and $\{(W,\mathcal{W}),(X,\mathcal{X}),u,P\}$ a RSCC with contraction. The Markov chain associated to the RSCC is regular, if and only if, there exists a point $\tilde{w} \in W$ such that

$$\lim_{n \to \infty} d(\sigma_n(\tilde{w}), w) = 0, \tag{2.14}$$

for any $w \in W$, where $\sigma_n(w) = \text{supp } Q^n(w, \cdot)$, where supp $\mu$ denotes the support of the measure $\mu$.

**Lemma 2.14** For any $m, n \in \mathbf{N}$, $w \in W$, we have

$$\sigma_{m+n}(w) = \overline{\bigcup_{w' \in \sigma_m(w)} \sigma_n(w')}, \tag{2.15}$$

where the line designates the topological aderence.

**Definition 2.15** Let $\{(W,\mathcal{W}),(X,\mathcal{X}),u,P\}$ be a RSCC. The RSCC is called *uniformly ergodic* if for any $r \in \mathbf{N}_+$ there exists a probability $P_r^{\infty}$ on $\mathcal{X}^r$ such that $\lim_{n \to \infty} \varepsilon_n = 0$, where

$$\varepsilon_n = \sup_{\substack{w \in W, r \in \mathbf{N}_+ \\ A \in \mathcal{X}^r}} |P_r^n(w, A) - P_r^{\infty}(A)|.$$

**Theorem 2.16** Let $(W,d)$ be a compact space. If the RSCC $\{(W,\mathcal{W}),(X,\mathcal{X}),u,P\}$ with contraction has regular associated Markov chain, then it is uniform ergodic.

## 3. THE GAUSS–KUZMIN TYPE THEOREM

**Proposition 3.1** The function $P(x,i) = P_i(x)$ from (1.10) defines a transition probability function from $(I, \mathcal{B}_I)$ to $(\mathbf{N}, \mathcal{P}(\mathbf{N}))$.

**Proof.** We have to verify that $\sum_{i \in \mathbf{N}} P(x,i) = 1$ for all $x \in I$. Since

$$P(x,i) = (x+1)\left(\frac{1}{x+i} - \frac{1}{x+i+1}\right),$$

then

$$\sum_{i \in \mathbf{N}} P(x,i) = (x+1) \cdot \frac{1}{x+1} = 1.$$

**Definition 3.2** Proposition 3.1 and relations (1.8) and (1.9) allows us to consider the random system with complete connections $\{(W,\mathcal{W}),(X,\mathcal{X}),u,P\}$, with

$$W = I, \quad \mathcal{W} = \mathcal{B}_I, \quad X = \mathbf{N}_+, \quad \mathcal{X} = \mathcal{P}(\mathbf{N}_+),$$

and

$$u(x,i) = \frac{1}{x+i}, \quad P_i(x) = \frac{x+1}{(x+i)(x+i+1)}, \quad x \in I, \; i \in \mathbf{N}_+.$$

**Remark 3.3** For $x_0 = 0$ and $P_0 = \lambda$ (see Theorem 2.5), the sequences $(\xi_n)_{n \in \mathbf{N}_+}$ and $(\zeta_n)_{n \in \mathbf{N}}$ associated with this RSCC coincide with the sequences $(a_n)_{n \in \mathbf{N}_+}$ and $(s_n)_{n \in \mathbf{N}}$, $s_0 = 0$ defined in (1.5) and (1.9).

**Proposition 3.4** The RSCC from Definition 3.1 is a RSCC with contraction and its associated Markov operator $U$ is regular with respect to $L(I)$ (the collection of all Lipschitz functions).

**Proof.** We have to verify the conditions from Definition 2.12. We have

$$\frac{d}{dx} P(x,i) = \frac{i^2 - i - (x+1)^2}{(x+i)^2 (x+i+1)^2}$$

and

$$\frac{d}{dx} u(x,i) = -\frac{1}{(x+i)^2}.$$

Hence, for any $x \in I$ and $i \in \mathbf{N}_+$, we have

$$\sup_{x \in I} \left| \frac{d}{dx} P(x,i) \right| < \frac{1}{i^2},$$

and

$$\sup_{x \in I} \left| \frac{d}{dx} u(x,i) \right| < \frac{1}{i^2}.$$

Thus, $R_1 < \infty$ and $r_1 < \infty$. To proof the regularity of $U$ with respect to $L(I)$, let us define recursively $x_{n+1} = \frac{1}{x_n + 2}$, $n \in \mathbf{N}_+$, with $x_0 = x$. Clearly, $x_{n+1} \in \sigma_1(x_n)$ and therefore Lemma 2.14 and an induction argument lead to the conclusion that $x_n \in \sigma_n(x_n)$, $n \in \mathbf{N}_+$. But $\lim_{n \to \infty} x_n = \sqrt{2} - 1$ for any $x \in I$. Hence

$$d(\sigma_n(x), \sqrt{2} - 1) \leq |x_n - \sqrt{2} - 1| \to 0, \; n \to \infty.$$

Finally, the regularity of $U$ with respect to $L(I)$ follows from Theorem 2.13.

Now, by the virtue of Theorem 2.16, the RSCC from Definition 3.1 is uniform ergodic. Moreover, $Q^n(\cdot, \cdot)$ converges uniformly to a probability measure $Q^\infty$ and that there exists two positive constants $q < 1$ and $k$ such that

$$\|U^n f - U^\infty f\|_L \leq k q^n \|f\|_L, \quad n \in \mathbf{N}_+, \; f \in L(I), \tag{3.1}$$

where

$$U^n f(\cdot) = \int_I f(y) Q^n(\cdot, dy), \quad U^\infty f = \int_I f(y) Q^\infty(dy), \tag{3.2}$$

and $\|\cdot\|_L$ is the norm over $L(I)$

$$\|f\|_L = \sup_{x \in I} |f(x)| + \sup_{x' \neq x''} \frac{|f(x') - f(x'')|}{|x' - x''|}.$$

**Proposition 3.5** The probability $Q^\infty$ coincide with the Gauss' measure $\gamma$ defined in (1.11).

**Proof.** By the virtue of uniqueness of $Q^\infty$ we have to prove that

$$\int_I Q(x, A) \gamma(dx) = \gamma(A), \quad A \in \mathcal{B}_I. \tag{3.3}$$

Since the intervals $[0,u) \subset [0,1]$ generate $\mathcal{B}_I$, it is suffices to check the above equation just for $A = [0, u)$, $0 < u \le 1$. We have

$$Q(x, [0,u)) = \sum_{\{i: u(x,i) \in A\}} P(x, i), \quad x \in I, \quad A \in \mathcal{B}_I.$$

Then

$$Q(x, [0, u)) = \sum_{i \ge \left\lfloor \frac{1}{u} - x \right\rfloor + 1} P(x, i) = \frac{x+1}{x + \left\lfloor \frac{1}{u} - x \right\rfloor + 1},$$

thus

$$\int_0^1 Q(x, [0, u)) \gamma(dx) = \frac{1}{\log 2} \int_0^1 \frac{dx}{x + \left\lfloor \frac{1}{u} - x \right\rfloor + 1} = \frac{1}{\log 2} \left( \int_0^{u^{-1} - \lfloor u^{-1} \rfloor} \frac{dx}{x + \left\lfloor \frac{1}{u} \right\rfloor + 1} + \int_{u^{-1} - \lfloor u^{-1} \rfloor}^1 \frac{dx}{x + \left\lfloor \frac{1}{u} \right\rfloor} \right)$$

$$= \frac{\log(1+u)}{\log 2} = \gamma([0,u)).$$

**Proposition 3.6** Let $\mu$ be an arbitrary non–atomic probability measure on $\mathcal{B}_I$. If

$$F_n(x) = F_n(x, \mu) = \mu(\tau^n < x), \quad x \in I, \quad n \in \mathbf{N}, \tag{3.4}$$

with $F_0(x) = \mu([0, x))$, then for any $n \in \mathbf{N}$, $F_n$ satisfies the following Gauss – Kuzmin type equation

$$F_{n+1}(x) = \sum_{i \in \mathbf{N}_+} \left( F_n\left(\frac{1}{i}\right) - F_n\left(\frac{1}{x+i}\right) \right), \quad x \in I. \tag{3.5}$$

**Proof.** Since $\tau^n = \dfrac{1}{a_{n+1} + \tau^{n+1}}$, it follows that

$$F_n(x) = \mu(\tau^n < x) = \sum_{i \in \mathbf{N}_+} \mu\left( \frac{1}{a_{n+1} + x} < \tau^n < \frac{1}{a_{n+1}} \right) = \sum_{i \in \mathbf{N}_+} \left( F_n\left(\frac{1}{i}\right) - F_n\left(\frac{1}{x+i}\right) \right).$$

Assuming that for some $m \in \mathbf{N}$, the derivative $F'_m$ exists everywhere in $I$ and is bounded, it is easy to see by induction that $F'_{m+n}$ exists and is bounded for all $n \in \mathbf{N}_+$. This allows us to differentiate (3.5) term by term, obtaining

$$F'_{n+1} = \sum_{i \in \mathbf{N}} \frac{1}{(x+i)^2} F'_n\left(\frac{1}{x+i}\right). \tag{3.6}$$

Further, write $f_n(x) = (x+1) F'_n$, $x \in I$, $n \in \mathbf{N}$. Then (3.6) becomes $f_{n+1} = U f_n$, $n \ge m$, with $U$ being the linear operator defined as

$$U f(x) = \sum_{i \in \mathbf{N}_+} P_i(x) f\left(\frac{1}{x+i}\right). \tag{3.7}$$

Now, let $\mu$ be a probability measure on $\mathcal{B}_I$ such that $\mu \ll \lambda$. Then it can be shown that

$$\mu(\tau^{-n}(A)) = \int_A \frac{U^n f_0(x)}{x+1} dx, \quad n \in \mathbf{N}, \quad A \in \mathcal{B}_I, \tag{3.8}$$

where

$$f_0(x) = (x+1) F'_0(x), \quad x \in I, \tag{3.9}$$

with $F'_0 = \dfrac{d\mu}{d\lambda}$.

**Theorem 3.7 (Gauss – Kuzmin Theorem)** Let $\mu$ be a probability measure on $\mathcal{B}_I$ such that $\mu \ll \lambda$. If the density $F'_0$ of $\mu$ is a Riemann integrable function, then

$$\lim_{n\to\infty}\mu(\tau^n < x) = \frac{1}{\log 2}\log(x+1),\ x \in I. \tag{3.10}$$

If the density $F_0'$ of $\mu$ is a Lipschitz function, then there exist two positive constants $q < 1$ and $k$ such that for all $x \in I$ and $n \in \mathbf{N}_+$

$$\mu(\tau^n < x) = \frac{1}{\log 2}(1+\theta q^n)\log(x+1), \tag{3.11}$$

where $\theta = \theta(\mu, n, x)$, with $|\theta| \leq k$.

**Proof.** Let $F_0'$ be a Lipschitz function. Then $f_0 \in L(I)$ and by the virtue of (3.2)

$$U^\infty f_0 = \int_I f_0(x) Q^\infty(dx) = \int_0^1 f_0(x)\gamma(dx) = \int_0^1 F_0'(x)dx = \frac{1}{\log 2}. \tag{3.12}$$

According to (3.1) there exist two constants $q < 1$ and $k$ such that

$$U^n f_0 = U^\infty f_0 + T^n f_0,\ n \in \mathbf{N}_+, \tag{3.13}$$

with $\|T^n f_0\|_L \leq kq^n$. Further, consider $C(I)$ the metric space of real continuous functions defined on $[0,1]$ with the supremum norm. Since $L(I)$ is a dense subset of $C(I)$ we have

$$\lim_{n\to\infty}|T^n f_0| = 0, \tag{3.14}$$

for $f_0 \in C(I)$. Therefore (3.14) is valid for measurable $f_0$ which is $\gamma$-almost surely continuous, that is for Riemann integrable $f_0$. Thus

$$\lim_{n\to\infty}\mu(\tau^n < x) = \lim_{n\to\infty}\int_0^x U^n f_0(u)\frac{1}{u+1}du = \int_0^x U^\infty f_0(u)\frac{1}{u+1}du =$$

$$= \frac{1}{\log 2}\log(u+1)\Big|_0^x = \frac{1}{\log 2}\log(x+1).$$

Equation (3.11) is equivalent with

$$F_n(x) = (1+\theta q^n)F_\infty(x),\ n \in \mathbf{N}_+,$$

which results from

$$U^n f_0(x) = (1+\theta q^n)U^\infty f_0(x),\ n \in \mathbf{N}_+,\ x \in I.$$

## 4. REFERENCES

[1] IOSIFESCU, M., *A very simple proof of a generalization of the Gauss-Kuzmin-Lèvy theorem on continued fractions, and related questions*, Rev. Roumaine Math. Pures Appl.37, pp.901 – 914, 1992
[2] IOSIFESCU, M., GRIGORESCU, S., *Dependence with Complete Connections and Its Applications*, Cambridge University Press, Cambridge, 1990
[3] IOSIFESCU, M., KRAAIKAMP, C., *Metrical theory of continued fractions*, Kluwer Academic Publishers, 2002
[4] KALPAZIDOU, S., *On a problem of Gauss-Kuzmin type for continued fraction with odd partial quotients*, Pacific J. Math. Vol.123, No.1, pp.103-114, 1986
[5] ROCKETT, A.M., SZÜSZ, P., *Continued fractions*, World Scientific, Singapore, 1992
[6] SCHWEIGER, F., *Ergodic theory of fibred systems and metric number theory*, Clarendon Press, Oxford, 1995
[7] SEBE, G.I., *On a Gauss-Kuzmin-type problem for a new continued fraction expansion with explicit invariant measure*, Proc. of the 3-rd Int. Coll "Math. in Engg. and Numerical Physics" (MENP-3), 7-9 October 2004 Bucharest, Romania, BSG Proceedings 12, Geometry Balkan Press, pp. 252-258, 2005